\documentclass[a4paper,12pt]{article}
\begin{document}

\title{Approximation and Schauder bases in M\"untz spaces $M_{\Lambda ,C}$ of continuous
functions.}
\author{S.V. Ludkowski}
\date{15 August 2016}
\maketitle

\begin{abstract}
M\"untz spaces $M_{\Lambda ,C}$ of continuous functions supplied
with the absolute maximum norm are considered in this article.
Fourier series approximations of functions in M\"untz spaces
$M_{\Lambda ,C}$ are studied. An existence of Schauder bases in
M\"untz spaces $M_{\Lambda ,C}$ is investigated. \footnote{key words
and phrases: Banach space; M\"untz space;
isomorphism; Schauder basis; Fourier series; approximation.  
\par Department of Applied Mathematics,
\par Moscow State Technological University MIREA,
\par av. Vernadsky 78, Moscow 119454, Russia
\par Ludkowski@mirea.ru}
\end{abstract}

\section{Introduction.}
The area of mathematics devoted to topological and geometric
properties of topological vector spaces is very important in
functional analysis (see, for example,
\cite{jarchb,koetheb,ludksmj,ludjms09,naribeckb}). Particularly,
studies of bases in Banach spaces compose a great part of it (see,
for example,
\cite{jarchb,koetheb,lindliorb,lusky92}-\cite{lusky04,schauder27,semadb,woytb}
and references therein). Many open problems remain for concrete
classes of Banach spaces. \par Among them M\"untz spaces $M_{\Lambda
,C}$ play very important role and there also remain unsolved
problems (see \cite{clarkerd43,gurlusb,schwartzb59} and references
therein). They are defined as completions of the linear span over
the real field $\bf R$ or the complex field $\bf C$ of monomials
$t^{\lambda }$ with $\lambda \in \Lambda $ on the segment $[0,1]$
relative to the absolute maximum norm, where $\Lambda \subset [0,
\infty )$, $~t\in [0,1]$. It was K. Weierstrass who in 1885 had
proved his theorem about polynomial approximations of continuous
functions on the segment. But the space of continuous functions
possesses the algebra structure. Later on in 1914 C. M\"untz had
considered generalizations so that his spaces generally had not such
algebraic structure.
\par A problem was whether they have bases \cite{gurizv66,schauder27}. Then a progress
was for M\"untz spaces satisfying the lacunary condition
$\underline{\lim}_{n\to \infty } \lambda _{n+1}/\lambda _n>1$ with
countable $\Lambda $, but in general this problem was unsolved
\cite{gurlusb,schwartzb59}. It is worth to mention that the
monomials $t^{\lambda }$ with $\lambda \in \Lambda $ generally do
not form a Schauder basis in the M\"untz space $M_{\Lambda ,C}$.
\par In this article results of investigations of the author on this problem are
presented.
\par In section 2 a Fourier approximation in M\"untz spaces
$M_{\Lambda ,C}$ of continuous functions on the unit segment $[0,1]$
supplied with the absolute maximum norm is studied. For this purpose
auxiliary Lemmas 3, 4 and Theorem 5 are proved. They are utilized
for reducing a consideration to a subclass of M\"untz spaces
$M_{\Lambda ,C}$ such that a set $\Lambda $ is contained in the set
of natural numbers $\bf N$ up to an isomorphism of Banach spaces. It
is proved that for M\"untz spaces satisfying the M\"untz condition
and the gap condition their functions belong to Weil-Nagy's class.
Then the theorem about existence of Schauder bases in M\"untz spaces
$M_{\Lambda ,C}$ under the M\"untz condition and the gap condition
is proved.
\par All main results of this paper are obtained for the first time.
They can be used for further studies of function approximations and
geometry of Banach spaces. This is important not only for progress
of mathematical analysis and functional analysis, but also in
different applications including measure theory and stochastic
processes in Banach spaces.

\section{M\"untz spaces $M_{\Lambda ,C}$.}
\par To avoid misunderstanding we first present our notation and
definitions.
\par {\bf 1. Notation.} Let $C([a,b],{\bf F})$ denote the Banach space
of all continuous functions $f: [a,b]\to {\bf F}$ supplied with the
absolute maximum norm \par $ \| f \|_{C[a,b]} := \max \{ |f(x)|:
x\in [a,b] \} $, \\ where $- \infty < a < b < \infty $, while ${\bf
F}=\bf R$ is the real field or ${\bf F}=\bf C$ the complex field.
\par Then $L_p((a,b), {\bf F})$ denotes the Banach space of
all Lebesgue measurable functions $f: (a,b) \to \bf F$ possessing a
norm as defined by the Lebesgue integral: $$ \| f \|_{L_p((a,b),
{\bf F})} := (\int_a^b |f(x)|^pdx)^{1/p} < \infty ,$$ where $1\le
p<\infty $ is a is a fixed number, and $- \infty \le a < b \le
\infty $.
\par Let $Q = (q_{n,k})$ be a lower triangular infinite matrix
with matrix elements
$q_{n,k}$ having values in the field ${\bf F}=\bf R$ or
${\bf F}=\bf C$ so that $q_{n,k}=0$
for each $k>n$, where $k, n$ are nonnegative integers.
To each $1$-periodic function $f: {\bf R}\to {\bf F}$
in the space $L_1([a,a+1],{\bf F})$ is counterposed a trigonometric polynomial
$$(1)\quad U_n(f,x,Q) := \frac{a_0}{2} q_{n,0} + \sum_{k=1}^n q_{n,k} (a_k
\cos (2\pi kx) + b_k \sin (2\pi kx)),$$
where $a_k=a_k(f)$ and $b=b_k(f)$ are the Fourier coefficients of a function
$f(x)$, whilst on $\bf R$ the Lebesgue measure is considered.
\par For measurable $1$-periodic functions $h$ and $g$ their
convolution is defined whenever it exists:
$$(2)\quad (h*g)(x) := 2 \int_{a}^{a+1}
h(x-t)g(t)dt.$$
Putting the kernel of the operator $U_n$ to be:
$$(3)\quad U_n(x,Q) := \frac{q_{n,0}}{2}+\sum_{k=1}^n q_{n,k} \cos (2\pi kx)$$
one gets
$$(4)\quad U_n(f,x,Q) = (f*U_n(,Q))(x)= (U_n(\cdot ,Q)*f)(x).$$
The norms of these operators are:
$$(5)\quad {\sf L}_n(Q) := \sup_{\| f \|_C = 1} \| U_n(f,x,Q) \|_C =  2
 \| U_n(x,Q) \|_{L_1} = 2 \int_a^{a+1} |U_n(t,Q)| dt,$$
where $\|* \|_C$ and $\| * \|_{L_1}$ denote norms on Banach spaces
$C([a,a+1],{\bf F})$ and $L_1([a,a+1],{\bf F})$ respectively, while
$a\in {\bf R}$ is a marked real number. These numbers ${\sf L}_n(Q)$
are called Lebesgue constants of a summation method.
\par Denote by $C_0([a,a+1],{\bf F})$ the subspace of continuous functions
\par $f: [a,a+1]\to \bf F$  satisfying the periodicity condition $f(0)=f(1)$.
\par As usually ${\bf N} = \{ 1, 2,... \} $ denotes the set of all natural numbers,
also $c_0({\bf F})$ stands for the Banach space of all converging to
zero sequences in $\bf F$ supplied with the absolute supremum norm.
\par Henceforth the Fourier summation methods given by $ \{ U_m : m \} $
which converge on $C_0([a,a+1],{\bf F})$
$$(6)\quad \lim_{m\to \infty } U_m(f,x,Q)=f(x)$$
uniformly in $x\in [a,a+1]$ will be considered.
\par {\bf 2. Definition.} Let $\Lambda $ be an increasing sequence in the set $(0,\infty )$.
\par The completion of the linear space containing all monomials
$a t^{\lambda }$ with $a\in \bf F$ and $\lambda \in \Lambda $
and $t\in [\alpha ,\beta ]$ relative to the absolute maximum norm
$$\| f(t) \|_{C[\alpha ,\beta ]} := \sup_{t\in [\alpha ,\beta ]} |f(t)|$$ is denoted by
$M_{\Lambda ,C}[\alpha ,\beta ]$, where $0\le \alpha <\beta <\infty
$. Particularly, for $[\alpha ,\beta ] =[0,1]$ it is also shortly
written as $M_{\Lambda ,C}$. We consider also its subspace
$$M^0_{\Lambda ,C}[0,1] := \{ f: f\in  M_{\Lambda ,C}[0,1]; ~
f(0)=f(1) \} $$ consisting of $1$-periodic functions.
\par Henceforth it is supposed that the sequence $\Lambda $ satisfies
the gap condition
\par $(2)\quad \inf_k \{ \lambda _{k+1} - \lambda _k \} =: \alpha _0>0$
and the M\"untz condition
\par $$(3)\quad \sum_{k=1}^{\infty } \frac{1}{\lambda _k}=:\alpha _1<\infty ,$$
where $\Lambda = \{ \lambda_k: k\in{\bf N} \} $.

\par Below Lemmas 3, 4 and Theorem 5 are proved about isomorphisms of M\"untz spaces
$M_{\Lambda ,C}$. Utilizing them we reduce our consideration to a
subclass of M\"untz spaces $M_{\Lambda ,C}$ such that a set $\Lambda
$ is contained in the set of natural numbers $\bf N$.
\par {\bf 3. Lemma.} {\it For each $0<\delta <1$ M\"untz spaces
$M_{\Lambda ,C}([0,1],{\bf F})$ and $M_{\Lambda ,C}([\delta ,1],{\bf
F})$ are isomorphic.}
\par {\bf Proof.} For every $0<\delta <1$ and $0<\epsilon \le 1$ and $f\in E$
the norms $\| f \|_{ C[0,1]}$ and $\epsilon \| f|_{[0, \delta ]}
\|_{ C[0,\delta ]} + \| f|_{[\delta , 1]} \|_{ C[\delta , 1]}$ are
equivalent. Due to the Remez-type inequality (see Theorem 6.2.2 in
\cite {borwerdb} and Theorem 7.4 in \cite{borerdjams1997}) for each
$\Lambda $ satisfying the M\"untz condition there is a constant
$\eta >0$ so that $\| h|_{[0, \delta ]} \|_{ C[0,\delta ]} \le \eta
\| h|_{[\delta , 1]} \|_{ C[\delta , 1]}$ for each $h\in M_{\Lambda
,C}$, where $\eta $ is independent of $h$. This implies that the
norms $\| h|_{[\delta , 1]} \|_{ C[\delta , 1]}$ and $ \| h \|_{
C[0, 1]}$ are equivalent on $M_{\Lambda ,C}[0,1]$.
 Each polynomial $a_1t^{\lambda _1}+...+a_nt^{\lambda _n}$
on $[\delta ,1]$ has the natural extension on $[0,1]$, where
$a_1,...,a_n\in {\bf F}$ are constants and $t$ is a variable. Thus
M\"untz spaces $M_{\Lambda ,C}[0,1]$ and $M_{\Lambda ,C}[\delta ,1]$
are isomorphic for each $0<\delta <1$.

\par {\bf 4. Lemma.} {\it M\"untz spaces $M_{\Lambda ,C}$ and
$M_{\Xi \cup (\alpha \Lambda +\beta ),C}$ are isomorphic for every $
\beta \ge 0$ and $\alpha >0$ and a finite subset $\Xi $ in
$(0,\infty )$.}
\par {\bf Proof.}  The set $\Lambda $ is infinite with $\lim_n \lambda _n = \infty $.
By virtue of Theorem 9.1.6 \cite{gurlusb} M\"untz space $M_{\Lambda
, C}$ contains a complemented isomorphic copy of $c_0({\bf F})$.
Therefore, $M_{\Lambda ,C}$ and $M_{\Xi \cup \Lambda ,C}$ are
isomorphic.
\par The isomorphism of
$M_{\alpha \Lambda ,C}$ with $M_{\Lambda ,C}$ follows from the
equality $$\sup_{t\in [0,1]} |f(t)| = \sup_{t\in [0,1]} |f(t^{\alpha
})|$$ for each continuous function $f: [0,1]\to \bf F$, since the
mapping $t\mapsto t^{\alpha }$ is the diffeomorphism of the segment
$[0,1]$ onto itself. Taking at first $\Lambda _1 = \Lambda \cup \{
\frac{\beta }{\alpha } \} $ and then $\alpha \Lambda _1$ we infer
that the M\"untz spaces $M_{\Lambda ,C}$ and $M_{\alpha \Lambda
+\beta ,C}$ are isomorphic.

\par {\bf 5. Theorem.} {\it Suppose that increasing sequences $\Lambda $ and
$\Upsilon $ of positive numbers satisfy the M\"untz condition and
the gap condition and $\lambda _n\le \mu _n$ for each $n$. If
\par $\sup_n (\mu _n-\lambda _n)=\delta $, where \par $\delta
<(8\sum_{n=1}^{\infty }\lambda _n^{-1})^{-1}$, \\ then $M_{\Lambda
,C}$ and $M_{\Upsilon, C}$ are the isomorphic Banach spaces.}
\par {\bf Proof.} Certainly the M\"untz spaces $M_{\Lambda ,C}$ and
$M_{\Upsilon, C}$ have isometric linear embeddings of
 into $M_{\Lambda \cup \Upsilon ,C}$. Consider a sequence of sets $\Upsilon _k$ satisfying
the following restrictions \par $(1)$ $\Upsilon _k= \{ \mu_{k,n}:
n=1,2,... \} \subset \Lambda \cup \Upsilon $ and $\mu _{k,n} \in \{
\lambda _n, \mu _n \} $ for each $k=0,1,2,...$ and $n=1,2,...$,
where $\Upsilon _0=\Lambda $;
\par $(2)$ $\mu _{k,n}\le \mu _{k+1,n}$ for each $k=0,1,2,...$ and
$n=1,2,...$;
\par $(3)$ $\{ \Delta _{k+1,n}: \Delta _{k+1,n}\ne 0; n=1,2,... \} $ is a monotone decreasing
subsequence converging to zero obtained from the sequence $\Delta
_{k+1,n} := \mu _{k+1,n}-\mu _{k,n}$ by elimination of zero terms;
\par $(4)$ $\{ m(k+1): k \} $ is a monotone increasing sequence with
$m(k+1) := \min \{ n: \mu _n -\mu _{k+1,n}\ne 0; \forall l<n ~ \mu
_l =\mu _{k+1,l} \} $.
\par For each $f\in M_{\Upsilon _k, C}$ we consider the power series
$f_1(t) = \sum_{n=1}^{\infty } a_n t^{\mu _{k+1,n}}$, where the
power series decomposition $f(t) = \sum_{n=1}^{\infty } a_n t^{\mu
_{k,n}}$ converges for each $0\le t <1$, since $f$ is analytic on
$[0,1)$ (see \cite{clarkerd43,gurlusb}). Then we infer that
$$f(t^2)-f_1(t^2) = \sum_{n=1}^{\infty } a_n t^{\mu
_{k,n}}u_n(t)\mbox{ with }u_n(t) := t^{\mu _{k,n}}- t^{\mu
_{k,n}+2\Delta _{k+1,n}}$$ so that $u_n(t)$ is a monotone decreasing
sequence by $n$ and hence \par $|f(t^2)-f_1(t^2)|\le 2
|u_{m(k+1)}(t)| |f(t)|$ \\ according to Dirichlet's criterium for
each $0\le t<1$. Therefore, the function $f_1(t)$ has the continuous
extension onto $[0,1]$ and $$ \| f-f_1 \| _{C([0,1],{\bf F})} \le 4
\| f \| _{C([0,1],{\bf F})} \Delta _{k+1,m(k+1)}/ \lambda
_{m(k+1)},$$ since the mapping $t\mapsto t^2$ is the order
preserving diffeomorphism of $[0,1]$ onto itself. Thus the series
$\sum_{l=1}^{\infty } a_n t^{\mu _{k+1,n}}$ converges on $[0,1)$.
Similarly to each $g_1\in M_{\Upsilon _{k+1}, C}$ a function $g\in
M_{\Upsilon _{k+1}, C}$ corresponds. \par This implies that there
exists the linear isomorphism $T_k$ of $M_{\Upsilon _k, C}$ with
$M_{\Upsilon _{k+1}, C}$ so that \par $ \| T_k -I \| \le 4 \Delta
_{k+1,m(k+1)}/ \lambda _{m(k+1)}$, $ ~ T_k: M_{\Upsilon _k, C}\to
M_{\Upsilon _{k+1}, C}$. \par Next we take the sequence of operators
\par $S_n := T_nT_{n-1}...T_0: M_{\Lambda , C}\to M_{\Upsilon _{n+1},
C}\subset M_{\Lambda \cup \Upsilon ,C}$. \par The space $M_{\Lambda
\cup \Upsilon ,C}$ is complete and the sequence $\{ S_n: n \} $
converges in the operator norm uniformity to an operator $S:
M_{\Lambda , C}\to M_{\Lambda \cup \Upsilon ,C}$ so that $ \| S-I \|
<1$, since \par $\sum_{k=0}^{\infty }\Delta _{k+1,m(k+1)}/ \lambda
_{m(k+1)}\le \delta \sum_{n=1}^{\infty }\lambda _n^{-1}<\infty ,$ \\
where $I$ denotes the unit operator. Therefore, the operator $S$ is
invertible. From Conditions $(1-4)$ it follows that $S(M_{\Lambda ,
C}) =M_{\Upsilon ,C}$.

\par {\bf 6. Remark.} In view of Lemmas 3, 4 and Theorem 5 it
suffices to consider a set $\Lambda $ satisfying the gap and M\"untz
conditions such that $\Lambda \subset \bf N$ up to an isomorphism of
M\"untz spaces. \par Next we recall necessary definitions and
notations of the Fourier approximation. Then auxiliary Proposition
10 is given which is used for proving Theorem 11 about the property
that for M\"untz spaces satisfying the M\"untz and gap conditions
their functions belong to Weil-Nagy's class.

\par {\bf 7. Notation.} Suppose that $(\psi (k): ~ k\in {\bf N} )$
is a sequence of non-zero numbers tending to zero, $\beta $ is a
marked real number. By $\sf F$ is denoted the set of all pairs
$(\psi , \beta )$, for which
$$(1)\quad {\cal D}_{\psi , \beta }(x) := \sum_{k=1}^{\infty }
\psi (k) \cos (2\pi kx+\beta \pi /2)$$ is the Fourier series of some
function belonging to $L_1[0, 1]$. Then ${\sf F}_1$ denotes the
family of all positive sequences $(\psi (k): ~ k\in {\bf N} )$
tending to zero with $\Delta_2\psi (k) := \psi (k-1) - 2\psi (k)
+\psi (k+1)\ge 0$ for each $k$ so that the series $$(2)\quad
\sum_{k=1}^{\infty } \frac{\psi (k)}{k} <\infty $$ converges.
\par An approximation of a function $f$ by the Fourier series
$S(f,x)$ is estimated by the function $$(3)\quad \rho _n (f,x) :=
f(x) - S_{n-1}(f,x),$$  where
$$(4)\quad S_n(f,x) := \frac{a_0}{2} + \sum_{k=1}^n (a_k
\cos (2\pi kx) + b_k \sin (2\pi kx))$$ is the partial Fourier sum.
That is a trigonometric polynomial approximating a summable (i.e.
Lebesgue integrable) $1$-periodic function $f\in L_1[0,1]$.

\par {\bf 8. Definitions.} Let $f\in L_1[a,a+1]$ and $S(f)$ be its Fourier series
with coefficients $a_k$ and $b_k$, let also $\psi (k)$ be an
arbitrary sequence real or complex. If the function $$f^{\psi
}_{\beta } := \sum_{k=1}^{\infty } [a_k(f)\cos (2\pi kx+\beta \pi
/2) + b_k(f) \sin (2\pi kx + \beta \pi /2)]/\psi (k)$$ belongs to
the space $L_1[a,a+1]$, then $ f^{\psi }_{\beta }= D^{\psi }_{\beta
} f$ is called the Weil ($(\psi ,\beta )$) derivative of $f$. \par
Let for a Banach space $\cal N$ of some functions on $[a,a+1]$:
$$C^{\psi }_{\beta }{\cal N}[a,a+1] := \{ f\in {\cal N}: ~
\exists ~ f^{\psi }_{\beta } \in C_0[a,a+1] \} ,$$ where $(\psi (k):
~ k)$ is a sequence with non-zero elements for each $k$ and  $\beta
$ is a real parameter.
\par In particular, let $C^{\psi }_{\beta }M[a,a+1]$
(or $C^{\psi }_{\beta }[a,a+1]$ for short) be the space of all
continuous $1$-periodic functions $f$ having a continuous Weil
derivative $f^{\psi }_{\beta } $, $~f^{\psi }_{\beta }\in C_0[0,1]$
and considered relative to the absolute maximum norm and such that
$$(1)\quad \| f^{\psi }_{\beta } \|_C := \max \{ |f^{\psi }_{\beta }
(t)|: ~ t\in [0,1] \} <\infty .$$ Particularly, for $\psi (k)
=k^{-r}$ there is the Weil-Nagy class \par $W^r_{\beta
,C}=W^r_{\beta }C[a, a +1] := \{ f: f\in C_0[a, a +1], \exists
f^{\psi }_{\beta }\in C_0[a, a+1] \} $.
\par Then let ${\cal E}_n (X) := \sup \{ \| \rho _n(f;x) \|
_{C[a,a +1]}: f \in X \} $,
\par $E_n(f) := \inf \{ \| f- T_{n-1} \|_{C[a ,a +1]}: T_{n-1}\in
{\cal T}_{2n-1} \} $,
\par $E_n(X) := \sup \{ E_n(f): f \in X
\} $, \\ where $X$ is a subset in $C_0[a ,a +1]$,
$${\cal T}_{2n-1} := \{ T_{n-1}(x) = \frac{c_0}{2} +
\sum_{k=1}^{n-1} (c_k \cos (2\pi kx) + d_k \sin (2\pi kx)); ~ c_k,
d_k \in {\bf R} \} $$ denotes the family of all trigonometric
polynomials $T_{n-1}$ of degree not greater than $n-1$.

\par {\bf 9. Note.}  We remind the following
definition: the family of all Lebesgue measurable functions $f:
(a,b)\to \bf R$ satisfying the condition
$$\| f \|_{L_{s,w}(a,b)} :=  \sup_{y>0}(y^s\mu \{ t: ~ t\in (a,b), ~ |f(t)|\ge y \} )^{1/s} <\infty
$$ is called the weak $L_s$ space and denoted by $L_{s,w}(a,b)$, where $\mu $ notates the Lebesgue
measure on the real field $\bf R$, $~0<s<\infty $, $~(a,b)\subset
{\bf R}$ (see, for example, \S 9.5 \cite{edwardsb}, \S IX.4
\cite{reedsim}, \cite{steinb}).

\par {\bf 10. Proposition.} {\it Suppose that an increasing sequence $\Lambda = \{
\lambda _n : n \} $ of natural numbers satisfies the M\"untz
condition. If $f\in M_{\Lambda ,C}[0,1]$, then $df(x)/dx\in
L_{1,w}(0,1)$.}
\par {\bf Proof.} Let $f\in M_{\Lambda ,C}[0,1]$. In view of
\cite{clarkerd43} a function $f$ has the analytic extension on
$\dot{B}_1(0)$ and the series
\par $(1)$ $f(z)=\sum_{n=1}^{\infty } a_nz^{\lambda _n}$\\
converges on $\dot{B}_1(0)$, where $\dot{B}_r(x) := \{ y: y\in {\bf
C}, |y-x|<r \} $ denotes the open disk in $\bf C$ of radius $r>0$
with center at $x\in {\bf C}$, where $a_n\in {\bf F}$ is an
expansion coefficient for each $n\in \bf N$. This also can be
lightly seen from Theorems 6.2.2, 6.2.3 and Corollary 6.2.4
\cite{gurlusb}, the Abel theorem 3 and the Cauchy-Hadamard formula
of a power series convergence radius and Theorem 5 in Subsection 20,
Section 6, Chapter II in \cite{shabatb}, since $\Lambda \subset {\bf
N}$. That is, the function $f$ has a holomorphic univalent extension
from $[0,1)$ on $\dot{B}_1(0)$. \par Using the Riemann integral we
have that
$$f(x)-f(0)=\int_0^x f'(t)dt\mbox{ for each }0\le x< 1\mbox{
and}$$
$$(2)\quad \lim_{x\uparrow 1}\int_0^x f'(t)dt=f(1)-f(0)$$ due to
Newton-Leibnitz' formula (see, for example, \S II.2.6 in
\cite{kamynb}), since $f(x)$ is continuous on $[0,1]$.  \par By
virtue of the uniqueness theorem for holomorphic functions (see, for
example, II.2.22 in \cite{shabatb}), applied to the considered case,
if a nonconstant holomorphic function $g$ on $\dot{B}_1(0)$ has a
set $E(g)=\{ x: g(x)=0, 0\le x< 1 \} $ of zeros in $[0,1)$, then
either $E(g)$ is finite or infinite with the unique limit point $1$.
Then we consider a linear function $s(x)=\alpha +\beta x$ with real
coefficients $\alpha $ and $\beta $, put $u(x)=f(x)+s(x)$ and choose
$\alpha $ and $\beta $ so that $u(0)=u(1)=0$.
\par On the other hand, $\min \Lambda = \lambda _1>0$ and hence
the function $f$ is nonconstant. The case $du(x)/dx=const $ is
trivial. So there remains the variant when $du(x)/dx$ is
nonconstant. Denote by $x_n$ zeros in $[0,1)$ of $du(t)/dt$ of odd
order so that $x_{n+1}>x_n$ for each $n\in {\bf N}$. Therefore,
$$(3)\quad \int^{x_{n+1}}_{x_n}u'(t)dt
\int^{x_{n+2}}_{x_{n+1}}u'(\tau )d\tau <0$$ for each $n\in {\bf N}$
according to Theorem II.2.6.10 in \cite{kamynb}. If $\{ x_n: n \} $
is a finite set, then from Formulas $(1,2)$ it follows that $u'\in
L_1[0,1]$ and hence $u'\in L_{1,w}[0,1]$. \par Consider now the case
when the set $\{ x_n: n \} $ is infinite. We take the closed disc $
V= \{ u\in {\bf C}: |u-1/2|\le 1/2 \} $, hence $[0,1]\subset V$.
According to Cauchy's formula 21$(5)$ in \cite{shabatb}
$$(4)\quad f'(z)= \frac{1}{2\pi i} \int_{\omega } \frac{f(\xi
)}{(\xi -z)^2}d\xi $$ for each $z\in Int (V)$, where $\omega $ is a
rectifiable path encompassing once a point $z$ in the positive
direction so that $\omega \subset Int (V)$, where $Int (V) $ denotes
the interior of $V$. For $3/4<x<1$ a circle can be taken with center
at $x$ and of radius $0<r<1-x$ with $r\uparrow (1-x)$ whilst
$x\uparrow 1$. \par Using the homotopy theorem (see Subsection 17,
Section 5, Chapter II in \cite{shabatb}) one can take simply the
circle $ \omega =\omega _q=
\partial V_q = \{ u\in {\bf C}: |u-1/2|= q \} $ with $q<1/2$, $ ~ r+x-1/2<q$ and
$q\uparrow 1/2$. From the standard parametrization of the circle
$\omega _q$ and the path integral (see Example 1 and Formula $(1)$
in Subsection 15, Section 5, Chapter II in \cite{shabatb}), Theorem
1 in Section IV.1.2 in \cite{kamynb} about continuity of an integral
depending on a parameter, the Cauchy criterion or the Weierstrass
criterion about uniform convergence of an improper integral
depending on a parameter in Section IV.1.3 in \cite{kamynb} it
easily follows that there exists a limit in the right side of $(4)$
if $q\uparrow 1/2$ and $1/3\le |z|\le x$, since $V$ is compact and
$\max_{z\in V} |f(z)|=: G <\infty $ by the Weierstrass theorem in
Section I.4.3 in \cite{kamynb}, the function $f$ is continuous on
$V$, $f$ is holomorphic on $V\setminus \{ 1 \} $, $u'(x_n)=0$ for
each $n$, where the sequence $\{ x_n : n \} $ has the limit point
$1$ in the considered case, $3/4<x<1$, $~ r+x-1/2<q<1/2$.
\par Then from the Cauchy estimate (see, for example, Subsection 20,
Section 6, Chapter II in \cite{shabatb}) of the integral $(4)$ it
evidently follows that
$$(5)\quad |f'(x)|\le
G/(2\pi (1-x))$$ for each $3/4<x<1$, hence $f'(x)\in L_{1,w}(3/4,1)$
and consequently $u'(t)\in L_{1,w}(3/4,1)$. Therefore, from the
inequality $(5)$ we infer that $$\sup_{t>0}y\mu \{ t: ~ t\in [0,1),
~ |f'(t)|\ge y \} \le $$ $$\sup_{t>0}y\mu \{ t: ~ t\in [0,3/4), ~
|f'(t)|\ge y \} + \sup_{t>0}y\mu \{ t: ~ t\in [3/4,1), ~ |f'(t)|\ge
y \} <\infty ,$$ where $\mu $ denotes the Lebesgue measure on
$[0,1]$. The latter means that $df(x)/dx\in L_{1,w}(0,1)$.

\par {\bf 11. Theorem.} {\it Let an increasing sequence $\Lambda = \{ \lambda _n: ~ n \}
\subset {\bf N}$ of natural numbers satisfy the M\"untz condition
and let $f\in M_{\Lambda ,C}[0,1]$. Then for each $0<\gamma <1$
there exists $\beta =\beta (\gamma )\in \bf R$ so that $v\in
W^{\gamma}_{\beta ,C}[0,1]$, where $v(t)=f(t)+(f(0)-f(1))t$ for each
$t\in [0,1)$ and $v$ is $1$-periodic on $\bf R$.}
\par {\bf Proof.} We have that $f(0)=0$, since $\lambda _n\ge 1$
for each $n$. Therefore, we consider $v(t)=f(t)- f(1)t$ on $[0,1) :=
\{ t: 0\le t<1 \} $ and take its $1$-periodic extension $v$ on $\bf
R$.
\par  According to Proposition
1.7.2 \cite{stepanetsb} (or see \cite{zygmb}) a function $h$ belongs
to $W^{\gamma }_{\beta }L_{\infty }(l,l+1)$ if and only if there
exists a function $\phi = \phi_{h, \gamma ,\beta }$ which is
$1$-periodic on $\bf R$ and Lebesgue integrable on $[0,1]$ such that
$$(1)\quad h(x) = \frac{a_0(h)}{2} + (\phi * {\cal D}_{\psi ,\beta
})(x),$$ where $a_0(h)= 2 \int_0^1 h(t)dt$ (see \S \S 7 and 8).
\par We take a sequence $U_n(t,Q)$ given by Formula 1$(3)$ so that
$$\lim_m q_{m,k} =1 \mbox { for each }k \mbox{ and } \sup_m
{\sf L}_m(Q) < \infty \mbox{ and }  \sup_{m,k} |q_{m,k}| <\infty $$
and write for short $U_n(t)$ instead of $U_n(t,Q)$. Under these
conditions the limit exists
\par $(2)$ $\lim_n (v*U_n)(x) =v(x)$ \\ in $L_{\infty }(0,1)$ norm for each
$v\in L_{\infty }((0,1),{\bf F})$ according to Chapters 2 and 3 in
\cite{stepanetsb} (see also \cite{barib,zygmb}).
\par On the other hand, Formula I$(10.1)$ \cite{stepanetsb} provides
\par $(3)$ $S[(y^{\psi _1}_{\bar{\beta }_1})^{\psi _2/\psi _1}_{\bar{\beta
}_2 - \bar{\beta }_1}]=S[y^{\psi _2}_{\bar{\beta}_2}]$, \\ where
$S[y]$ is the Fourier series corresponding to a function $y\in
L^{\psi _2}_{\bar{\beta}_2}$, when $(\psi _1,\bar{\beta }_1)\le
(\psi _2,\bar{\beta }_2)$.
\par Put $\theta (k)=k^{\gamma -1}$ for all $k\in \bf N$. Then
${\cal D}_{\theta , - \beta }\in L_1(0,1)$ for each $\beta \in \bf
R$ due to Theorems II.13.7, V.1.5 and V.2.24 \cite{zygmb} (or see
\cite{barib}). This is also seen from chapters I and V in
\cite{stepanetsb} and Formulas $(1)$ and $(3)$ above. In view of
Dirichlet's theorem (see \S 430 in \cite{fichtb}) the function
${\cal D}_{\theta , -\beta }(x) $ is continuous on the segment
$[\delta ,1-\delta ]$ for each $0<\delta <1/4$.
\par According to formula 2.5.3.$(10)$ in \cite{prbrmarb1}
$$\int_0^{\infty } x^{\alpha -1} {\sin (bx)\choose{\cos (bx)}} dx =
b^{-\alpha }\Gamma (\alpha )
 {\sin (\pi \alpha /2) \choose{\cos (\pi \alpha /2)}}$$
 for each $b>0$ and $0<Re (\alpha )<1$. On the other hand,
the integration by parts gives: $$\int_a^{\infty } x^{\alpha -1}
{\sin (bx)\choose{\cos (bx)}} dx = b^{-1}a^{\alpha -1} {\cos
(ab)\choose{- \sin (ab)}} - b^{-1}(\alpha -1)\int_a^{\infty }
x^{\alpha -2} {- \cos (bx)\choose{\sin (bx)}} dx$$ for every $a>0$,
$b>0$ and $0<Re (\alpha )<1$. From formulas V$(2.1)$, theorems
V.2.22 and V.2.24 in \cite{zygmb} (see also \cite{bruijnb,olverb})
we infer the asymptotic expansions
$$\sum_{n=1}^{\infty } n^{-\alpha }\sin (2\pi nx) \approx  (2\pi x)^{\alpha -1} \Gamma (1-\alpha )
\cos (\pi \alpha /2) + \mu x^{\alpha },$$
$$\sum_{n=1}^{\infty } n^{-\alpha }\cos (2\pi nx) \approx  (2\pi x)^{\alpha -1} \Gamma (1-\alpha )
\sin (\pi \alpha /2)+ \nu x^{\alpha }$$ in a small neighborhood
$0<x<\delta $ of zero, where $0<\delta <1/4$, $0<\alpha <1$, $\mu $
and $\nu $ are real constants. Taking $\beta =\alpha = 1-\gamma $ we
get that ${\cal D}_{\theta , -\beta }(x) \in L_{\infty }(0,1).$
\par Remind that for Lebesgue measurable
functions $f: {\bf R}\to {\bf R}$ and $g: {\bf R}\to {\bf R}$ there
is the equality $\int_{-\infty }^{\infty } f(x-t)\chi _{[0,\infty
)}(x-t)g(t)\chi _{[0,\infty )}(t) dt = \int_0^x f(x-t)g(t)dt$ for
each $x>0$ whenever this integral exists, where $\chi _A$ denotes
the characteristic function of a subset $A$ in $\bf R$ such that
$\chi _A(y)=1$ for each $y\in A$, also $\chi _A(y)=0$ for each $y$
outside $A$, $y\in {\bf R}\setminus A$. In particular, if $0<x\le
T$, where $0<T<\infty $ is a constant, then $\int_0^x f(x-t)g(t)dt=
\int_0^{\infty } f(x-t)\chi _{[0,T]}(x-t)g(t)\chi _{[0,T]}(t)dt$
(see also \cite{fichtb,kolmfomb}). This can be applied to formula
1$(2)$ putting $\alpha =0$ there and with the help of the equality
\par $\int_0^1 f(x-t)g(t)dt = \int_0^xf(x-t)g(t)dt + \int_0^{1-x}
f_1((1-x)-v)g_1(v)dv$ \\ for each $0\le x\le 1$ and $1$-periodic
functions $f$ and $g$ and using also that $ \| f|_{[a,b]} \| \le \|
f|_{[0,1]} \| = \| f_1|_{[0,1]} \|$ for the considered here types of
norms for each $[a,b]\subset [0,1]$, where $f_1(t)=f(-t)$ and
$g_1(t)=g(-t)$ for each $t\in \bf R$, since \par $\int_x^1
f(x-t)g(t)dt = \int_0^{1-x} f(v-1+x) g(1-v)dv$.
\par Mention that according to the weak Young inequality
\par $(4)$ $ \| \xi *\eta  \|_r \le K_{p,q}
\| \xi \|_p \| \eta  \|_{q,w}$ \\
for each $\xi \in L_p$ and $\eta \in L_{q,w}$, where $1\le p, r \le
\infty $, $1\le q<\infty $ and $p^{-1} + q^{-1} =1 + r^{-1}$,
$~K_{p,q}>0$ is a constant independent of $\xi $ and $\eta $ (see
Theorem 9.5.1 in \cite{edwardsb}, \S IX.4 in \cite{reedsim}).
\par With the help of Proposition 10 and $(2)$ we define
the function $s(x)$ such that
$$s(x)=\lim_{\eta \downarrow 0}
\lim_n \eta ^{-1} \int_0^{\eta }(({\cal D}_{\theta , - \beta
}*U_n)*v')(x-t)dt.$$ In virtue of the weak Young inequality $(4)$
and Proposition 10 this function $s$ is in $L_{\infty }(0,1)$.
Therefore $\phi_{v, \gamma ,\beta }=s$ and $D^{\psi }_{\beta }v=s$
according to $(1)$ and $(3)$. Thus $v\in W^{\gamma }_{\beta
}L_{\infty }(0,1)$. On the other hand, $v$ is analytic on $(0,1)$,
$1$-periodic and continuous on $[0,1)$, consequently, $s$ is
analytic on $(0,1)$ and $1$-periodic. Together with Formulas $(1-3)$
the latter implies that $\| s \|_{L_{\infty }(0,1)} = \| s
\|_{C(0,1)}$ and hence $v\in W^{\gamma }_{\beta ,C } [0,1]$.

\par {\bf 12. Lemma.} {\it  If an increasing sequence $\Lambda $
of natural numbers satisfies the M\"untz condition, also $0<\gamma
<1$,
\par $X= \{ v: ~ v \mbox{ is } 1-\mbox{periodic and} ~ \forall t\in [0,1) ~ v(t) = f(t) +(f(0)- f(1))t,$\par $ f\in
M_{\Lambda ,C}; \| v \|_{C[0,1]}\le 1 \} $,
\par then a positive constant
$\omega  = \omega (\gamma )$ exists so that
$$(1)\quad E_n(X)\le {\cal E}_n(X) \le \omega n^{-\gamma }\ln n
$$ for each natural number $n\in \bf N$.}
\par {\bf Proof.} Due to Theorem 11 the inclusion is valid
$v\in W^{\gamma }_{\beta ,C}[0,1]$ for the $1$-periodic extension
$v$ of $v(t) = f(t)+(f(0)-f(1))t$ on $[0,1)$ for each $f\in
M_{\Lambda ,C}$.
\par Then estimate $(1)$ follows from Theorems 3.12.3 and 3.12.3' in \cite{stepanetsb}.

\par {\bf 13. Lemma.} {\it If $\psi \in {\sf F}_1$
and $(\psi , \beta )\in {\sf F}$ (see \S 7), then $C^{\psi
,0}_{\beta }[0,1] := \{ f\in C^{\psi }_{\beta }[0,1]:
\int_0^1f(t)dt=0 \} $ is the Banach space relative to the norm given
by the formula:
$$(1)\quad \| f \| _{C^{\psi }_{\beta }[0,1]} := \| f^{\psi }_{\beta }
\|_{C[0,1]} .$$}
\par {\bf Proof.} We have that $C^{\psi }_{\beta }[0,1]$ is the $\bf F$-linear
space and hence $C^{\psi ,0}_{\beta }[0,1]$ is such also as the kernel of
the linear functional $\phi (f) := \int_0^1f(t)dt$, since each $f\in C^{\psi }_{\beta }[0,1]$
is integrable. Therefore, the assertion of this lemma
follows from Propositions I.8.1 and I.8.3
\cite{stepanetsb}, since each $f\in C^{\psi }_{\beta }[0,1]$
has the convolution representation:
$$(2)\quad f(x) = \frac{a_0(f)}{2} + 2 \int_0^1 f^{\psi }_{\beta }(x-t)
{\cal D}_{\psi ,\beta }(t)dt$$ for each $x\in [0,1]$, but $a_0(f)=0$
for each $f\in C^{\psi ,0}_{\beta }[0,1]$, while the convolution
$h*u$ is continuous for each $h\in C[0,1]$ and $u\in L_1[0,1]$ so
that $ \| h*u \|_C \le 2 \| h \|_C \| u \|_{L_1}$, where ${\cal
D}_{\psi ,\beta }$ is defined by Formula 7$(1)$.
\par To prove Theorem 15 about existence of a Schauder basis
the following proposition is useful.

\par {\bf 14. Proposition.} {\it Let $X$ be a Banach space over $\bf R$
and let $Y$ be its Banach subspace so that they fulfill conditions
$(1-4)$ below:
\par $(1)$ there is a sequence $(e_i: i \in {\bf N})$ in $X$
such that $e_1,...,e_n$ are linearly independent vectors and $ \|
e_n \|_X =1$ for each $n$ and
\par $(2)$ there exists a Schauder basis $(z_n: n\in {\bf N})$
in $X$ such that \par $z_n = \sum_{k=1}^n b_{k,n}e_k$ for each $n\in
\bf N$, where $b_{k,n}$ are real coefficients;
\par $(3)$ for every $x\in Y$ and $n\in \bf N$ there exist
$x_1,...,x_n\in {\bf R}$ so that
\par $ \| x - \sum_{i=1}^n x_ie_i \|_X \le s(n) \| x \| $, \\
where $s(n)$ is a strictly monotone decreasing positive function
with \par $\lim_{n\to \infty } s(n)=0$ and
\par $(4)$ $u_n =\sum_{l=m(n)}^{k(n)} u_{n,l} e_l$, \\
where $u_{n,l}\in \bf R$ for each natural numbers $k$ and $l$, where
a sequence $(u_n: n \in {\bf N})$ of normalized vectors in $Y$ is
such that its real linear span is everywhere dense in $Y$ and $1\le
m(n)\le k(n)<\infty $ and $m(n)<m(n+1)$ for each $n\in \bf N$.
\par Then $Y$ has a Schauder basis.}
\par {\bf Proof.} Without loss of generality one can select and enumerate
\par $(5)$ vectors $u_1$,...,$u_n$ so that they are linearly independent in $Y$
for each natural number $n$. By virtue of Theorem $(8.4.8)$ in
\cite{naribeckb} their real linear span $span_{\bf R}
(u_1,....,u_n)$ is complemented in $Y$ for each $n\in \bf N$. Put
$L_{n,\infty } := cl_X span_{\bf R} (u_k: k\ge n)$ and $L_{n,m}:=
cl_X span_{\bf R} (u_k: n\le k\le m)$, where $cl_XA$ denotes the
closure of a subset $A$ in $X$, where $span_{\bf R}A$ denotes the
real linear span of $A$. Since $Y$ is a Banach space and $u_k\in Y$
for each $k$, then $L_{n,\infty } \subset Y$ and $L_{n,m}\subset Y$
for each natural numbers $n$ and $m$. Then we infer that \par
$L_{n,j} \subset span_{\bf R} (e_l: m(n)\le l \le k_{n,j} )$, where
$k_{n,j} := \max (k(l): n\le l \le j)$. \par Take arbitrary vectors
$f\in L_{1,j}$ and $g\in L_{j+1,q}$, where $1\le j<q$. Therefore,
there are real coefficients $f_i$ and $g_i$ such that
\par $ f=\sum_{i=1}^{k_{1,j}} f_ie_i $ and
 \par $g=\sum_{i=m(j+1)}^{k_{j+1,q}} g_ie_i$. Hence due to condition
 $(2)$:
\par $ \| f - \sum_{i=1}^{m(j)} f_ie_i \| _X \le s(m(j)) \| f \| $
and \par $ \| g - \sum_{i=k_{1,j}+1}^{k_{j+1,q}} g_ie_i \|_X \le
s(k_{1,j}+1) \| g \|_X$. \par On the other hand,
\par $ f=  \sum_{i=1}^{m(j)} f_ie_i + \sum_{i=m(j)+1}^{k_{1,j}} f_ie_i $, consequently,
\par $ \| f^{[j+1]} \| \le s(m(j+1))) \| f \| $, where
\par $f^{[j+1]} := \sum_{i=m(j+1)}^{k_{1,j}} f_ie_i$ and  $\sum_{i=a}^bf_ie_i := 0$, when
$a>b$.  \par When $0<\delta <1/4$ and $s(m(j)+1)< \delta $ we infer
using the triangle inequality that $ \| f^{[j+1]} - h \|_X \le
\delta \| f^{[j+1]} \|_X /(1-\delta ) \le \delta s(m(j+1)-1) \| f
\|_X /(1-\delta )$ for the best approximation $h$ of $f^{[j+1]}$ in
$L_{j+1,\infty }$, since $m(j)<m(j+1)$ for each $j$. Therefore, the
inequality $ \| f - g \|_X \ge \| f-f^{[j+1} \|_X -  \| f^{[j+1]} -
g \|_X$ and $s(n)\downarrow 0$ imply that there exists $n_0$ such
that the inclination of $L_{1,j}$ to $L_{j+1,\infty }$ is not less
than $1/2$ for each $j\ge n_0$. Condition $(4)$ implies that
$L_{1,n_0}$ is complemented in $Y$. By virtue of Theorem 1.2.3
\cite{gurlusb} a Schauder basis exists in $Y$.

\par {\bf 15. Theorem.} {\it If an increasing sequence $\Lambda $
of positive numbers satisfies the M\"untz condition and the gap
condition, then the M\"untz space $M_{\Lambda ,C}[0,1]$ has a
Schauder basis.}
\par {\bf Proof.} By virtue of Lemma 4 and Theorem 5 it is sufficient to prove an existence
of a Schauder basis in the M\"untz space $M_{\Lambda ,C}$ for
$\Lambda \subset {\bf N}$. According to \S 4 the Banach spaces
$M^0_{\Lambda ,C}$ and $M_{\Lambda ,C}$ are isomorphic.
\par The functional
$$(1)\quad \phi (h) := \int_0^1 h(\tau )d\tau $$ is continuous on $C^{\psi
}_{\beta }[0,1]$, where $\psi $ and $\beta $ satisfy conditions of
Lemma 13. Then $coker (\phi ) = {\bf F}$. Therefore, $C^{\psi
}_{\beta }[0,1] = {\bf F}\oplus C^{\psi ,0}_{\beta }[0,1]$.
 \par In view of Theorem 6.2.3 and Corollary 6.2.4
\cite{gurlusb} each function $g\in M_{\Lambda ,C}[0,1]$ has the
analytic extension on $\dot{B}_1(0)$ and hence
$$(2)\quad g(z) = \sum_{n=1}^{\infty } c_n z^{\lambda _n} =
\sum_{k=1}^{\infty } p_k u_k(z)$$ are the convergent series on the
unit open disk $\dot{B}_1(0)$ in $\bf C$ with center at zero (see \S
12), where $\Lambda \subset \bf N$ and $c_n=c_n(g)\in \bf N$,
$~p_n=p_n(g)=c_1+...+c_n$, $~ u_1(z) := z^{\lambda _1}$, $~
u_{n+1}(z) := z^{\lambda _{n+1}} - z^{\lambda _n}$ for each
$n=1,2,...$.
\par Take the finite dimensional subspace $X_n := span_{\bf
R}(u_1,...,u_n)$ in $X:=M^0_{\Lambda ,C} $, where $n\in \bf N$. Due
to Lemma 4 the Banach space $X\ominus X_n$ exists and is isomorphic
with $M_{\Lambda ,C}$. \par Consider the trigonometric polynomials
$U_m(f,x,Q)$ for $f\in X$, where $m=1, 2,...$ (see \S \S 1 and 8
also). Put $Y_{K,n}$ to be the $C[0,1]$ completion of the linear
span $span_{\bf R} (U_m(f,x,Q): (m,f)\in K)$, where $K\subset {\bf
N}\times (X\ominus X_n)$, $m\in \bf N$, $f\in (X\ominus X_n)$.
\par There exists a countable subset $\{ f_n: n\in {\bf N} \} $ in
$X$ such that $f_n = {\cal D}_{\psi ,\beta }*g_n$ with $g_n\in
L(0,1)$ for each $n\in \bf N$ and so that $span_{\bf R}\{ f_n: n\in
{\bf N} \}$ is dense in $X$, since $X$ is separable. From Formulas
$(1,2)$ and Theorem 11 and Lemma 12 we infer that a countable set
$K$ and a sufficiently large natural number $n_0$ exist so that the
Banach space $Y_{K,n_0}$ is isomorphic with $(X\ominus X_{n_0})$ and
$Y_{K,n_0}|_{(0,1)}\subset W^{\gamma }_{\beta ,C}(0,1)$, where
$0<\gamma <1$ and $\beta =1-\gamma $. Thus the Banach space
$Y_{K,n_0}$ is the $C[0,1]$ completion of the real linear span of a
countable family $(s_l: l\in {\bf N})$ of trigonometric polynomials
$s_l$.
\par Without loss of generality this family can be refined by
induction such that $s_l$ is linearly independent of
$s_1,...,s_{l-1}$ over $\bf F$ for each $l\in \bf N$. With the help
of transpositions in the sequence $\{ s_l: l \in {\bf N} \} $, the
normalization and the Gaussian exclusion algorithm we construct a
sequence $ \{ r_l: l\in {\bf N} \} $ of trigonometric polynomials
which are finite real linear combinations of the initial
trigonometric polynomials $\{ s_l : l\in {\bf N} \} $ and satisfying
the conditions
\par $(3)$ $ \| r_l \|_{C(0,1)} =1$ for each $l$;
\par $(4)$ the infinite matrix having $l$-th row of the form
$...,a_{l,k}, b_{l,k}, a_{l,k+1}, b_{l,k+1},...$ for each $l\in \bf
N$ is upper trapezoidal (step), where $$r_l(x) = \frac{a_{l,0}}{2} +
\sum_{k=m(l)}^{n(l)} [a_{l,k} \cos (2\pi kx) + b_{l,k} \sin (2\pi
kx) ]$$ with $a_{l,m(l)}^2+b_{l,m(l)}^2>0$ and
$a_{l,n(l)}^2+b_{l,n(l)}^2>0$, where $1\le m(l)\le n(l)$, $deg
(r_l)=n(l)$, or $r_1(x)= \frac{a_{1,0}}{2}$ when $deg (r_1)=0$;
$a_{l,k}, b_{l,k}\in {\bf R}$ for each $l\in {\bf N}$ and $0\le k\in
{\bf Z}$.
\par Then as $X$ and $Y$ in Proposition 14 we take $X=C[0,1]$ and
$Y=Y_{K,n_0}$. In view of Proposition 14 and Lemma 4 the Schauder
basis exists in $Y_{K,n_0}$ and hence also in $M_{\Lambda ,C}[0,1]$.

\end{document}